\documentclass[12pt,a4paper,oneside,onecolumn]{article}
\title{On Pursell - Shanks type results}
\author{Lecomte P.B.A. and Zihindula Mushengezi E.}
\usepackage{etex}
\usepackage[]{latexsym}
\usepackage[T1]{fontenc}
\usepackage[latin1]{inputenc}
\usepackage[english]{babel}
\usepackage{graphics}
\usepackage[]{amssymb}
\usepackage{multicol}
\usepackage{amsmath}
\usepackage{graphicx}
\usepackage{amsfonts}
\usepackage[all]{xy}
\usepackage{pifont}
\usepackage{microtype}
\usepackage{answers}
\usepackage{tabularx}
\usepackage{float}
\usepackage{color}

\usepackage{listings}
\definecolor{darkWhite}{rgb}{0.94,0.94,0.94}
\usepackage{palatino}

\newcount\exenum     \newdimen\numindent
\newcount\qqnum     \newdimen\qqmarge
\newcount\sqnum     \newdimen\sqmarge

\outer\def\bye{%
\vskip 0pt\@endmulticol\@endgroup              
\ifanswer \let\next=\exobye@                   
\else     \let\next=\@exobye                   
\fi\next}

\usepackage[Glenn]{fncychap}

\usepackage{pstricks}
\usepackage{pgf,tikz}
\usepackage{pstricks-add}
\usepackage{pst-plot}
\usetikzlibrary{arrows}
\usepackage{pgfplots}
\usepackage{pgfkeys}

\definecolor{qqqqqq}{rgb}{0,0,0}

\definecolor{xdxdff}{rgb}{0.49,0.49,1}
\definecolor{qqwuqq}{rgb}{0,0.39,0}
\definecolor{qqqqff}{rgb}{0,0,1}
\definecolor{ttttff}{rgb}{0.2,0.2,1}
\definecolor{uququq}{rgb}{0.25,0.25,0.25}

\renewcommand{\epsilon}{ \varepsilon}

\newcommand{\euro}{\texteuro{}}

\newcommand{\pre}{{\bf Proof.\ }}

\newcommand{\cl }{\mathcal }

\newtheorem{theo}{Theorem}[section]
\newtheorem{prop}[theo]{Proposition}

\newcounter{exercice}

\definecolor{fbase}{rgb}{0.8,0.8,1}
\definecolor{fgris}{gray}{0.6}
\definecolor{frouge}{HTML}{DC143C}
\definecolor{fvert}{rgb}{0.6,1,0.6}
\definecolor{fbleu}{rgb}{0.4,0.4,1}
\definecolor{fjaune}{HTML}{DCDC14}
{\endMakeFramed}
\makeatletter
\DeclareRobustCommand\sfrac[1]{\@ifnextchar/{\@sfrac{#1}}%
                                            {\@sfrac{#1}/}}
\def\@sfrac#1/#2{\leavevmode\kern.1em\raise.5ex
         \hbox{$\m@th{\fontsize\sf@size\z@\selectfont#1}$}
         \kern-.1em/\kern-.15em\lower.55ex
          \hbox{$\m@th{\fontsize\sf@size\z@\selectfont#2}$}}

\DeclareRobustCommand{\Efrac}[2]{{\displaystyle\begingroup
\raise2ex\hbox{$\m@th{#1}$}\endgroup\@@over \lower1ex
\hbox{$\m@th{#2}$}}}
\makeatother

\numberwithin{equation}{section}

\newtheorem{Exc}{Exercice}
\Newassociation{correction}{Soln}{mycor}
\Newassociation{indication}{Indi}{myind}

\def\exo#1{\futurelet\testchar\MaybeOptArgmyexoo}
\def\MaybeOptArgmyexoo{\ifx[\testchar \let\next\OptArgmyexoo
                        \else \let\next\NoOptArgmyexoo \fi \next}
\def\OptArgmyexoo[#1]{\begin{Exc}[#1]\normalfont}
\def\NoOptArgmyexoo{\begin{Exc}\normalfont}

\newcommand{\finexo}{\end{Exc}}
\newcommand{\flag}[1]{}

\tikzset{
xmin/.store in=\xmin, xmin/.default=-3, xmin=-3,
xmax/.store in=\xmax, xmax/.default=3, xmax=3,
ymin/.store in=\ymin, ymin/.default=-3, ymin=-3,
ymax/.store in=\ymax, ymax/.default=3, ymax=3,
}

\newcommand{\entete}[1]

\begin{document}
 \Opensolutionfile{mycor}[ficcorex]
 \Opensolutionfile{myind}[ficind]
 \entete{\'Enoncés}


\maketitle

\begin{abstract}
We prove a Lie-algebraic characterization of vector bundle for the Lie algebra $\mathcal{D}(E,M),$ seen as ${\rm C}^\infty(M)-$module, of all linear operators acting on sections of a vector bundle $E\to M$. 

We obtain similar result for its Lie subalgebra $\mathcal{D}^1(E,M)$ of all linear first-order differential operators.

Thanks to a well-chosen filtration, $\mathcal{D}(E,M)$ becomes  $\mathcal{P}(E,M)$ and we prove that  $\mathcal{P}^1(E,M)$ characterizes the vector bundle without the hypothesis of being seen as ${\rm C}^\infty(M)-$module. 

We prove that the Lie algebra $\mathcal{S}(\mathcal{P}(E,M))$ of symbols of linear operators acting on smooth sections of a vector bundle $E\to M,$ characterizes it. 
To obtain this, we assume that $\mathcal{S}(\mathcal{P}(E,M))$ is seen as ${\rm C}^\infty(M)-$module.

We obtain a similar result with the Lie algebra $\mathcal{S}^1(\mathcal{P}(E,M))$ of symbols of first-order linear operators without the hypothesis of being seen as a ${\rm C}^\infty(M)-$module.
\end{abstract}

\section{Introduction}

A Lie algebra $ \cl D (M) $ characterizes the manifold $ M $ when for any manifold $ N, $ the Lie algebras $ \cl D (M) $ and $ \cl D (N) $ are isomorphic if, and only if, the manifolds $ M $ and $ N $ are diffeomorphic.\\
We then usually speak of a Pursell-Shanks type result since the characterization of manifolds by Lie algebras dates back to the result of Pursell and Shanks obtained in 1954. They used the Lie algebra of all $\rm C^\infty$-vector fields with compact support on a manifold $M.$ See \cite{PurShan}. \\
A whole series of papers followed using sub-Lie algebras of $Vect(M),$ the space of all $\rm C^\infty$-vector fields of $M.$ 
Note that object of interest is often Lie algebras of all infinitesimal automorphisms for several geometric structures on $M$.
Without being exhaustive, we can cite
\begin{enumerate}
\item In 1974, Koriyama \cite{Kor} obtained a Pursell-Shanks type result with considering a submanifold as a geometric structure.
\item Ameniya \cite{Ame}, in 1975, proved the corresponding result in the case of complex structures.
\item Omori \cite{Omo} (1976), studied the case of volume structures, symplectic structures, contact structures and fibering structures with compact fibers and obtained results of Lie-algebraic characterization of manifold.
\item Fukui (1981), see \cite{Fuk}, proved the corresponding result in case of Lie algebras of $G$-invariant $\rm C^\infty$-vector fields with compact support on paracompact, connected, free $G$-manifolds when $G$ is a compact connected semi-simple Lie group such that the automorphism group of its Lie algebra is connected.
\item Fukui and Tomita, see \cite{FukTom}, proved in 1982 a result corresponding to the case of foliated structures. In particular, they obtained a Pursell-Shanks type result for two special cases of one codimension foliation.
\item Rybicki obtained in 1990, see \cite{Ryb1}, a Pursell-Shanks type result for an another special case of one codimension foliation. In 1993, \cite{Ryb2}, he proved similar result for Riemannian foliated structures.
\end{enumerate}
Pursell and Shanks initial result has been generalized several times.\\
There are at least two ways of generalizing Pursell and Shanks' result. A first way consists not in characterizing a manifold but rather a vector bundle. We can cite, in this regard, Lecomte's paper, see \cite{Lec1} 1981, in which he proved that a vector bundle is characterized by the Lie algebra of all its infinitesimal automorphisms.\\
We can also generalize Pursell and Shanks result by considering Lie algebras larger than that of vector fields of a manifold.
Grabowski and Poncin followed this way by developing methods and concepts, see \cite{GraPon1,GraPon4}, allowing to conclude to the characterization of manifolds by $\cl D(M)$, the Lie algebra of all the differential operators on $M$ and by that of their symbols.\\
Our purpose in this work is to generalize the result of Lecompte and those of Grabowski and Poncin. This is also a way of generalizing the result of Pursell and Shanks.\\
In \cite{GraPon4}, Grabowski and Poncin tried to obtain a such generalization with the Lie algebra $\cl D(E,M)$ of all linear operators acting on smooth sections of a one rank vector bundle $E\to M.$ \\
They proved that $\cl D(M)$ and $\cl D(E,M)$ are isomorphic Lie algebras.\\
This implies that for a vector bundle $E\to M$ of rank $n=1,$ we cannot obtain any Pursell-Shanks type result with the algebra $\cl D(E,M).$\\
The powerfull methods developed by Grabowski and Poncin work for quantum Poisson algebras, see  \cite{GraPon2}, for instance. But $\cl D(E,M)$ is a quantum Poisson algebra if and only if $E\to M$ is a linear bundle.\\
For vector bundles of rank $n>1,$ we need therefore to define a more general concept and the methods that go with it. This is why we propose in the following section the notion of "quasi quantum Poisson algebra". 

\section{Quasi quantum Poisson algebra}

We start by building a theoretical framework allowing to study the Lie algebra $\cl D(E,M).$ \\ 
Let $\cl D=\bigcup_{i\geq 0} \cl D^i$ be an associative filtered algebra with unit over a commutative field $K$ of characteristic 0. We put $\cl D^0=\cl A$ and extend the filtration on $\mathbb{Z},$ by setting $\cl D^i=\{0\},$ for $i<0.$
\begin{enumerate}
\item We say the Lie algebra $\cl D$ with its commutators bracket is a  \textit{quasi quantum Poisson} algebra if
  \[
      \cl D^{i}\subset \cl D^{i+1},\quad \cl D^{i} \cl D^{j} =\cl D^{i+j} \text{ and }\quad [\cl D^{i} ,\cl D^{j}]\subset \cl D^{i+j}.
  \] 
In this case, $\cl A$ is an associative sub-algebra  of $\cl D,$ and also a Lie sub-algebra of $\cl D.$ 
The associative $\mathbb{R}-$algebra $\cl A$ is then called the \textit{basis algebra} or simply the \textit{basis} of the quasi quantum Poisson algebra $\cl D.$
\item The algebra $\cl D$ is said non-singular if there exists a Lie sub-algebra $\cl D^1_{\cl A}$ of  $\cl D^1$ such that  the center of $\cl A$ is $Z(\cl A)=[\cl D^1_{\cl A},Z(\cl A)].$
\item $\cl D$ is \textit{symplectic} if $Z(\cl D),$  the center of $\cl D,$ is reduced to constants. The space of constants being identified to $K$ by $k\in K\mapsto k\cdot 1\in\cl D,$ where $1$ designates the unit of the $\mathbb{R}-$algebra $\cl D.$
\item The \textit{centralizer} of $ad(Z(\cl A))$ in $Hom_K(\cl D,\cl D)$, denoted by $\cl C(\cl D),$ is defined by
  \[
    \Psi\in\cl C(\cl D)\Leftrightarrow \Psi([T,u])=[\Psi(T),u], \forall u\in Z(\cl A), \forall T\in\cl D
  \]  
\item $\cl D$ is \textit{quasi-distinguishing} if
 \[
   [T,u]=0 ,\forall u\in Z(\cl A)\Rightarrow T\in \cl A 
\]
and if for any integer $i,$
  \[
   \{T\in\cl D:[T,Z(\cl A)]\subset \cl D^i\}=\cl D^{i+1}.
  \]
\item By definition, $T\in Nil(\cl D)\Leftrightarrow \forall D\in\cl D, \exists n\in\mathbb{N}: (ad T)^n(D)=0.$     
\end{enumerate}
\begin{prop}\label{centralisateur}
 Let $\cl D$  be a quasi quantum Poisson algebra, non-singular and quasi-distinguishing. Then any $\Psi\in\cl C(\cl D)$ respects the filtration and we have 
  \[
    \Psi(u)=\Psi(1)u
  \]
for all $u\in Z(\cl A).$ 
  \end{prop} 
\pre
Consider $\Psi\in\cl C(\cl D).$ Then for all $A\in\cl D^0$ and all $u\in Z(\cl A),$ we have
\[
  [\Psi(A),u]=\Psi([A,u])=0 
\]
since $[A,u]=0.$ As $\cl D$ is quasi-distinguishing, it comes that $\Psi(A)\in\cl A$. Assume now that $\Psi(\cl D^{i})\subset \cl D^i. $ Then, for any $T\in\cl D^{i+1}, u\in Z(\cl A),$ we have
 \[
  [\Psi(T),u]=\Psi([T,u])\in\Psi(\cl D^i)\subset \cl D^i. 
 \]
Therefore, $\cl D$ being quasi-distinguishing, $\Psi(T)\in\cl D^{i+1}.$\\
Note first that for $T\in\cl D^1,$ the map $A\longmapsto \widehat{T}(A)=[T,A]$ is not usually a derivation of $\cl D^0,$ but  it verifies the relation 
\begin{equation}\label{(*)}
 \widehat{T}(AB)=\widehat{T}(A)B +A \widehat{T}(B)  
\end{equation}
since in general, $ad\Delta$ is a derivation of the associative structure of the  quasi quantum Poisson algebra $\cl D,$ for any $\Delta\in \cl D.$
Hence, for any $\Psi\in\cl C(\cl D), T\in\cl D^1, u\in Z(\cl A),$ we have, on the one hand
  \begin{eqnarray*}
    \Psi([T,u^2]) & = & \Psi(\widehat{T}(u^2))\\
                  & = & \psi(\widehat{T}(u)u+u\widehat{T}(u))\\
                  & = & 2\Psi(u\widehat{T}(u)).
  \end{eqnarray*}
The last equality comes because $T\in \cl D^1$ induces $\widehat{T}(u)\in \cl D^0.$  
On the other hand, by using the fact that $\Psi(T)\in\cl D^1,$  we have
  \begin{eqnarray*}
     \Psi([T,u^2]) & = & [\Psi(T),u^2]
                    =  2u\Psi(\widehat{T}(u)).
  \end{eqnarray*}
Therefore, 
  \begin{equation}\label{(**)}
    \Psi(u\widehat{T}(u))= u \Psi(\widehat{T}(u)),
  \end{equation}
for any $T\in\cl D^1, u\in Z(\cl A).$\\
Next, observe that for any $v\in Z(\cl A),$ the restrictions of $\widehat{\widehat{T}(v)T}$ and $\widehat{T}(v)\widehat{T}$ on $Z(\cl A)$ coincide. Indeed, for any $w\in Z(\cl A),$
   \begin{eqnarray*}
\widehat{\widehat{T}(v)T}(w)=
                          [\widehat{T}(v)T,w] & = & \widehat{T}(v)Tw -w \widehat{T}(v)T\\
                                              & = & \widehat{T}(v)Tw - \widehat{T}(v)w T\\
                                              & = & \widehat{T}(v)\widehat{T}(w).
   \end{eqnarray*}
In general, for any element $A$ in $\cl D^0,$ $A\widehat{T}$ does not verify the previous relation (\ref{(*)}) but $A\widehat{T}$ and $\widehat{AT}$ coincide on $Z(\cl A).$\\
Therefore, by replacing $T$ with $AT, A\in \cl D^0$ and $u$ with $u+w,w\in Z(\cl A), $ the relation (\ref{(**)}) becomes
  \[
    \Psi((u+w)\widehat{AT}(u+w))=(u+w)\Psi(\widehat{AT}(u+w)).
  \]
This is equivalent to
  \[
  \Psi(uA\widehat{T}(w))+\Psi(wA\widehat{T}(u))=u\Psi(A\widehat{T}(w))+w \Psi(A\widehat{T}(u)).
  \]
For $A=\widehat{T}(w)$ the last equality becomes
  \[
 \Psi(u(\widehat{T}(w))^2)+\Psi(w\widehat{T}(w)\widehat{T}(u))
 =
 u\Psi((\widehat{T}(w))^2)+w \Psi(\widehat{T}(w)\widehat{T}(u))
  \]
This relation becomes for $T\in\cl D^1_{\cl A},$
  \begin{equation}\label{(***)}
   \Psi(u(\widehat{T}(w))^2)=u\Psi((\widehat{T}(w))^2)   
  \end{equation}  
Indeed, if $T\in\cl D^1_{\cl A}$ then
   \begin{eqnarray*}
 \Psi(w\widehat{T}(w)\widehat{T}(u))& = & \Psi(w\widehat{T}(u)\widehat{T}(w)) \\
                                    & = & \Psi(w\widehat{\widehat{T}(u)T}(w)) \\
                                    & = & w\Psi(\widehat{\widehat{T}(u)T}(w))
                                      =   w \Psi(\widehat{T}(w)\widehat{T}(u)),
   \end{eqnarray*}
the penultimate equality resulting from (\ref{(**)}).\\
   Let us consider now the following set
  \[
 \cl J= \{ v\in Z(\cl A): \Psi(vu)=u\Psi(v), \forall u\in Z(\cl A)\}.
  \]
The sub-vector space  $\cl J$ is an ideal of $Z(\cl A).$ Indeed, for any  $u\in Z(\cl A), v\in\cl J$ we have $uv\in\cl J$ since for any $w\in Z(\cl A)$,
  \[
    w\Psi(uv) = (uw)\Psi(v)=\Psi((uw)v)=\Psi((uv)w).
  \] 
The previous relation (\ref{(***)}) shows that 
$Rad(\cl J)\supset[\cl D^1_{\cl A},Z(\cl A)].$ \\
The algebra $\cl D$ being non-singular, we conclude that 
     $Rad(\cl J)=Z(\cl A)$ and thus $\cl J=Z(\cl A).$
Hence the relation 
  \[
    \Psi(u)=u\Psi(1),
  \]
for any $\Psi\in\cl C(\cl D), u\in Z(\cl D^0).$\hfill $\blacksquare$

\begin{prop}\label{iso modules alg quasi filtré}
      Let $\cl D_1$ and $\cl D_2$ two quasi quantum  Poisson algebras, non-singular and quasi-distinguishing. Then, for any isomorphism $\Phi: \cl D_1\longrightarrow \cl D_2$ of Lie algebras such that $\Phi(Z(\cl A_1))=Z(\cl A_2)$ respects the filtration, its restriction to $\cl A_1$ is an isomorphism  of Lie algebras and its restriction to $Z(\cl A_1)$ is of the form
      \[
       \Phi|_{Z(\cl A_1)}=\kappa \Psi,
      \]
where $\Psi:Z(\cl A_1)\to Z(\cl A_2)$ is an isomorphism of associative algebras.
\end{prop}
\pre
Let us first observe that for all $ A \in \cl A_1, $ the endomorphisms of vector space $ \gamma_A: T \longmapsto AT $ and $ \Phi \circ \gamma_A \circ \Phi^{- 1}, $ of $ \cl D_1 $ and $ \cl D_2 $ respectively, belong to $ \cl C (\cl D_1) $ and $ \cl C (\cl D_2) $ as appropriate. Indeed, for any $A \in\cl A_1, T\in\cl D_1$ and $u\in Z(\cl A_1)$
  \[
   \gamma_A([T,u]) =  A Tu-Au T
                     =  [\gamma_A(T),u];
  \]
and for any $v\in Z(\cl A_2),$ since $\Phi^{-1}(v)\in Z(\cl A_1),$ we have for  $D\in\cl D_2,$ 
\[
  \Phi\circ\gamma_A\circ\Phi^{-1}([D,u])= 
                               [\Phi\circ\gamma_A\circ\Phi^{-1}(D),u].
\]
Therefore, for any $w\in Z(\cl A_2),$ in virtue of the previous Proposition \ref{centralisateur}, we obtain
  \begin{equation}\label{(iv)}
    \Phi\circ\gamma_A\circ\Phi^{-1}(w)=\Phi(A\Phi^{-1}(1))w\in \cl A_2 
  \end{equation}
We deduce that for any $A\in\cl A_1,$
  \begin{equation}\label{mboxv}
     \Phi(A\Phi^{-1}(1))\in \cl A_2 
  \end{equation}
Setting $v=\Phi^{-1}(w), A=u$ and $\lambda=\Phi^{-1}(1)$ in the previous relation (\ref{(iv)}), we obtain 
  \begin{eqnarray*}
   \Phi(uv) & = & \Phi(u\lambda)\Phi(v)\\ 
            & = & \Phi(\lambda u)\Phi(v)=\Phi(\lambda^2)\Phi( u)\Phi(v). 
  \end{eqnarray*}
In particular, for $u=1$ and $v=\lambda,$ the last equality above shows that the element $\Phi(\lambda^2),$ which belongs to the center of $\cl D_2,$ is invertible in $\cl A_2.$ 
Therefore, putting $\kappa^{-1}=\Phi(\lambda^2),$ we have that
 $\Psi:Z(\cl A_1)\to Z(\cl A_2): u\mapsto \kappa^{-1}\Phi(u),$  is an isomorphism of associative algebras. \\
Observe that 
 \begin{eqnarray*}
   \Psi(\lambda\Psi^{-1}(\kappa)) & = & \Psi(\lambda)\kappa\\ 
                                  & = & \kappa^{-1}\Phi(\lambda )\kappa=1.
  \end{eqnarray*}
Thus $\lambda=\Phi^{-1}(1)$ is invertible in $\cl A_1$ and we can thus deduce from the relation $(\ref{mboxv})$ the inclusion $\Phi(\cl A_1)\subset\cl A_2.$\\ 
Let us assume now that $\Phi(\cl D^i_{1})\subset\cl D^i_{2}$ for $i\in\mathbb{N}.$ Let $T\in\cl D^{i+1}_{1}.$ \\
We have
 \[
  [\Phi(T),Z(\cl A_2)]= \Phi([T,Z(\cl A_1)])\subset \cl D^{i}_{2}
 \]
Hence $\Phi(T)\in\cl D^{i+1}_{2}$ and then $\Phi$ respects the filtration. \hfill $\blacksquare$

\section{Differential operators of a vector bundle}

In the \cite{LecLeuZih}, we have studied Lie algebras made up of differential operators acting on smooth functions belonging to $ {\rm C}^\infty (E,\mathbb{R}),$ for a given vector bundle $ E \to M. $ In the following lines, we are interested in differential operators acting on the sections of a vector bundle.

Let $ E \to M $ be a vector bundle and denote by $ \Gamma (E) $ the space of its smooth sections. The space $ \Gamma (E) $ being a $ {\rm C}^\infty (M) - $ module, let us set
  \[
  \gamma_u:\Gamma(E)\to\Gamma(E):s\mapsto us, \forall u\in{\rm C}^\infty(M).
  \]
We have that $ \gamma_u $ is an endomorphism of the space $ \Gamma (E). $  
Let us then set
  \[
    \cl A(E,M):=\cl D^0(E,M)=\{T\in End(\Gamma(E)):[T,\gamma_u]=0,\,\forall u\in{\rm C}^\infty(M)\}
  \]
and for any integer $k\geq1,$
  \[  
    \cl D^k(E,M)= \{T\in End(\Gamma(E))| \forall u\in{\rm C}^\infty(M): [T,\gamma_u]\in    
    \cl D^{k-1}(E,M)\};
  \]
 where $[\cdot,\cdot]$ is the commutators bracket of $End(\Gamma(E)).$ \\
We can then state the following result.
\begin{prop}\label{localité ds dem}
 The elements of the  $\mathbb{R}-$space $\cl D^k(E,M),\, k\in\mathbb{N},$ are differential operators  of $\leq k$-order on $\Gamma(E)$. 
\end{prop}
\pre
The proof is done by induction on $ k $ and can be found in many papers.
The initial condition is well verified for $ k = 0 $ and suppose by induction hypothesis that the result is true for the elements of $\cl D^i(E,M), i\leq k-1.$ Let $T\in\cl D^k(E,M)$ and a section $s\in\Gamma(E)$ such that $j^1_x(s)=0$. We then have
  \[
        T(s)(x)=0,
  \]
  using the fact that the section $ s $ can then be decomposed into 
  \[
    s=\sum_{i\leq r} f_{i_1}\cdots f_{i_{k+1}}s_i, \quad f_{i_j}\in{\rm C}^\infty(M), s_i\in\Gamma(E),
  \]
 where the functions $f_{i_j}$ vanish at $x.$
  What completes the proof.\hfill $\blacksquare$\\

We can now state the following result where we group together the first properties of the spaces we have just defined. Let us specify before some notations.
In the following, the algebra  of endomorphisms fields of $ E, $ is denoted by $ gl(E). $ We thus have
\[
  \Gamma(Hom(E,E))=gl(E).
\]
The subalgebra of $ gl (E) $ of endomorphisms of $ E $ having a zero trace is denoted by $ sl (E) $ and it is the  derived ideal of the Lie algebra $ gl (E). $ We know that the center of $ gl (E) $
is given by
 \[
   Z(gl(E))={\rm C}^\infty(M)\, id
 \]
\begin{prop}\label{proprietes dem}
 The following relations are verified
\begin{enumerate}
  \item[(1)] $\cl A(E,M)=\cl D^0(E,M)=gl(E)$
  \item[(2)] $ \cl D^{i}(E,M)\subset \cl D^{i+1}(E,M),\, \cl D^{i}(E,M)\cdot \cl D^{j} (E,M) \subset\cl D^{i+j}(E,M)$
  \item[(3)] $ [\cl D^{i}(E,M) ,\cl D^{j}(E,M)]\subset \cl D^{i+j}(E,M).$
\end{enumerate}  
\end{prop}    
\pre
The relation $ (1) $ comes from the definition of the space $ \cl D^0 (E, M) $ and from the preceding Proposition \ref{localité ds dem}.\\
The two inclusions of $ (2) $ are proved by induction. For the first, we have that the initial condition is direct for $ i = 0. $ 
Let us suppose by induction hypothesis that this inclusion is true for $i=k\in\mathbb{N}$ and consider $T\in\cl D^{k+1}(E,M).$
We then have, for any function $u\in{\rm C}^\infty(M),$
\[
    [T,\gamma_u]\in\cl D^{k+1}(E,M),
\]
by the induction hypothesis.
This shows that $T\in\cl D^{k+2}(E,M)$ and the first inclusion is well established. 
For the second inclusion of $ (2), $ the induction is done on the sum $k=i+j.$
Remember that we have equality  
\begin{equation}\label{(v)}
   [T\circ T',T'']=T\circ[T',T'']+[T,T'']\circ T'
\end{equation}
for any $T\in\cl D^{p}(E,M),$ $T'\in\cl D^{q}(E,M)$ and $T''\in\cl D^{r}(E,M)$ with $p,q,r\in\mathbb{N}.$ Particularly, for $T,T'\in\cl D^{0}(E,M)$ and $T''=\gamma_u,$ where $u\in{\rm C}^\infty(M),$ the relation (\ref{(v)}) allows to obtain the initial condition for $k=0.$ \\ 
Assume by induction hypothesis that the inclusion is true for $i+j<k$ and consider $T\in\cl D^{p}(E,M)$ and $T'\in\cl D^{q}(E,M)$ with $p+q=k.$
We then have, by applying the induction hypothesis, that
 \[
  T\circ[ T',\gamma_u,]\,\, ,\,\, [T,\gamma_u]\circ T'\in\cl D^{p+q-1}(E,M),
 \]
for any $u\in{\rm C}^\infty(M).$
Hence, by virtue of the previous (\ref{(v)}) relation
 \[
    T\circ T'\in\cl D^{p+q}(E,M).
 \]
 The relation $ (3) $ is a consequence of the second inclusion of $ (2). $\hfill $\blacksquare$\\
Consider the Lie algebra defined by
  \[
   \cl D(E,M) =\bigcup_{k\geq 0} \cl D^k(E,M).
  \]
It is therefore a quasi quantum Poisson algebra whose elements will be called differential operators of the vector bundle $E\to M.$\\
\begin{prop}\label{alg quasi P q non-sing,sympl, quasi-q}
   Let $E\to M$ a  vector bundle of rank $n>1.$ The quasi quantum Poisson algebra $\cl D(E,M)$ is non-singular, quasi distinguishing and symplectic. 
\end{prop}
\pre
For the non-singularity, observe for a covariant derivation  $\nabla$  of the vector bundle $E$ and a vector field $X$ on $M,$ we have, for any $u\in{\rm C}^\infty(M),$
 \[
   [\nabla_X,\gamma_u]=\gamma_{X(u)}=\gamma_{[X,u]}.
 \]
Observe that this last bracket is the one defined in $ \cl D (M), $ the Lie algebra of the differential operators of $ M. $ The non-singularity of $ \cl D (E, M) $ then comes from that of $ \cl D (M), $
by setting 
 \[
   \cl D^1_{\cl A}(E,M)=\{T\in\cl D^1(E,M): [T,Z(\cl A(E,M))]\subset Z(\cl A(E,M))\}.
 \] 
Indeed, thus defined, the space $ \cl D^1_{\cl A}(E,M)$ is a Lie subalgebra of $ \cl D^1 (E, M) $ by the Jacobi identity.\\
If in a domain of trivialization $U\in M$, $T\in\cl D^k(E,M)$ is written
  \[
  T=\sum_{|\beta|<k}B_{\beta}\partial^\beta+\sum_{|\alpha|=k}A_\alpha\partial^\alpha,
  \] 
with $T_\alpha,T_\beta\in{\rm C}^\infty(U,gl(n,\mathbb{R})),$ then the relation $[T,\gamma_u]=0$ induces
   \[
   \sum_{|\alpha|=k}[A_\alpha\partial^\alpha,\gamma_u]=0.
   \]
Consequently, by considering the terms of maximum order, we obtain for any $u\in{\rm C}^\infty(M),$
 \[
   \sum_{|\alpha|=k} A_\alpha\circ\gamma_{\partial_{i}(u)}\partial^{\alpha_i}=0, 
 \] 
with for $i\in [1,m],$  $\alpha_i=\alpha-e_i.$ We deduce that $A_\alpha\circ\gamma_{\partial_{i}(u)}=0,$ and then $A_\alpha=0.$
Hence the following implication
 \[
   ([T,\gamma_u]=0 ,\forall u\in Z(\cl A(E,M)))\Rightarrow T\in \cl A(E,M) ;
 \]
which allows to conclude that $\cl D(E,M)$ is quasi distinguishing.\\
 We also have the following inclusion 
  \[
  Z(\cl D(E,M))\subset \{T\in\cl D(E,M):[T,Z(\cl A(E,M))]=0\}.
  \]
Therefore, $\cl D(E,M)$  being quasi-distinguishing, we have 
 \[
   Z(\cl D(E,M))\subset Z(\cl A(E,M)).
 \]    

In addition, for $u\in{\rm C}^\infty(M),$ the following equality 
\[
[\cl D^1_{\cl A}(E,M),\gamma_u]=0
\] 
implies that $u$ is constant. \hfill $\blacksquare$\\

We will now state results of Lie-algebraic characterization of vector bundles. The following theorem is taken from \cite{Lec3}. 

\begin{theo}\label{ende carct e}
Let $ E \mapsto M $ and $ F \mapsto M $ be two vector bundles of respective ranks $n,n'>1$ with $H^1(M,\mathbb{Z}/2)=0$. The Lie algebras $ gl (E) $ and $ gl (F) $ (resp. $ sl (E) $ and $ sl (F) $) are isomorphic if and only if the vector bundles $ E $ and $ F $ are isomorphic.
\end{theo}

\begin{theo}\label{dem carct e}
Let $ E \mapsto M $ and $ F \mapsto M $ be two vector bundles of respective ranks $n,n'>1$ with $H^1(M,\mathbb{Z}/2)=0$. 
The Lie algebras $ \cl D (E, M) $ and $ \cl D (F, M), $ seen as $ {\rm C} ^\infty (M) - $ modules, are isomorphic if and only if the vector bundles $ E $ and $ F $ are isomorphic.
\end{theo}
\pre
Let $\Phi: \cl D(E,M)\to \cl D(F,N)$ be an isomorphism of Lie algebras and of ${\rm C}^\infty(M)-$modules. Since these quasi quantum Poisson algebras are symplectic, we have 
\[
  \Phi({\rm C}^\infty(M))={\rm C}^\infty(N),
\]
by identifying $u\in{\rm C}^\infty(X)$ with $\gamma_u$ for $X=M$ or $X=N.$ Indeed, for any $u\in{\rm C}^\infty(M),$ we have
 \begin{eqnarray*}
 \Phi(\gamma_u) & = & \Phi(u\gamma_1)\\
                & = & u\Phi(\gamma_1).
 \end{eqnarray*}
The algebras in question here being non-singular and quasi-distinguishing, we can therefore apply the previous Proposition \ref{iso modules alg quasi filtré} to obtain the following equality 
  \[
    \Phi(gl(E))=gl(F).
  \]
The previous Theorem \ref{ende carct e} allows then to conclude. \hfill $\blacksquare$

\section{The quantum Poisson algebra $\cl P(E,M)$}

Let $ E \to M $ be a vector bundle of rank $ n> 1. $ In this section, we propose another filtration for the quasi quantum Poisson algebra $ \cl D (E, M), $ so as to make quantum Poisson algebra. 
We pose 
\[
  \cl P(E,M)=\cup_{k\geq 0}\cl P^k(E,M),
\]
with, by definition,
\begin{enumerate}
\item[-]  $ \cl P^0(E,M)=\{\gamma_u: u\in{\rm C}^\infty(M)\}$
\item[-]  $\cl P^{k+1}(E,M)=\{T\in End(\Gamma(E))|\forall u\in{\rm C}^\infty(M): [T,\gamma_u]\in\cl P^k(E,M)\}.$
 \end{enumerate}
We now give some properties of the algebra $ \cl P (E, M). $ Unlike $ \cl D (E, M), $ the filtration of $ \cl P (E, M) $ turns out to be that of a quantum Poisson algebra.
\begin{prop}
For any integer $j,k\in\mathbb{N},$ we have
\begin{enumerate}
  \item $\cl P^k(E,M)\subset \cl P^{k+1}(E,M)$ and $\cl P^{j}(E,M)\cdot\cl P^{k}(E,M)\subset\cl P^{j+k}(E,M)$
  \item $[\cl P^{j}(E,M),\cl P^{k}(E,M)]\subset\cl P^{j+k-1}(E,M)$
\end{enumerate}
\end{prop}
 \pre
The relation $ (1) $ can be proved exactly as for $ \cl D (E, M). $ We prove the last inclusion. Let us begin with $ \cl P^i (E, M) = \{0 \} $ for $ i <0 $ and do a recurrence on $ j + k. $ The result being true for $ j + k = 0, $ suppose the same is true for $j+k<p.$ Let then $D\in \cl P^j(E,M), T\in\cl P^k(E,M)$ with $j+k=p.$ We have for any $u\in{\rm C}^\infty(M),$ 
\[
[[D,\gamma_u],T]+[D,[T,\gamma_u]]\in\cl P^{j+k-2}(E,M).
\]
This comes by induction hypothesis, since $[D,\gamma_u]\in\cl P^{j-1}(E,M)$ and $[T,\gamma_u]\in\cl P^{k-1}(E,M),$ by definition. 
Therefore, we have
 \[
  [[D,T],\gamma_u]\in\cl P^{j+k-2}(E,M).
 \] 
 And this completes the proof of the proposition. \hfill $\blacksquare$\\
The link between the Lie algebras $ \cl D (E, M) $ and $ \cl P (E, M) $ is given in the following statement.
\begin{prop} Let $E\to M$ be a vector bundle.
We have, for any $k\in \mathbb{N},$
\begin{enumerate}
  \item $\cl P^k(E,M)\subset\cl D^k(E,M)\subset\cl P^{k+1}(E,M).$
  \item $\cl P(E,M)=\cl D(E,M)$
\end{enumerate}
\end{prop}
\pre
Observe that the relation $ (2) $ is a corollary of the relation $ (1). $ We thus prove $ (1), $ in the lines which follow and that is sufficient. The first inclusion being obvious, let's establish the second by induction. We have for every $ A \in \cl D^0 (E, M) $ and every function $u\in{\rm C}^\infty(M),$
\[
  [A,\gamma_u]\in\cl P^0(E,M);
\]  
and the inclusion is established for $k=0.$\\
Assume that it is also true for $k\in\mathbb{N}.$ Let $T\in \cl D^{k+1}(E,M).$ We then have
 \[
   [T,\gamma_u]\in\cl D^k(E,M)\subset\cl P^{k+1}(E,M)
 \]
for any  $u\in{\rm C}^\infty(M),$ this last inclusion being true by induction hypothesis. \\
We deduce the following inclusion
  \[
    \cl D^{k+1}(E,M)\subset\cl P^{k+2}(E,M).
  \]
And this concludes the proof of the proposition. \hfill  $\blacksquare$ 

\begin{prop}\label{pe non sing symple pas quasi dist}
The quantum Poisson algebra $ \cl P (E, M) $ is non-singular and symplectic but it is not quasi-distinguishing. We have the following results
  \begin{equation}\label{(v*)}
    [T,\gamma_u]=0,\forall u\in{\rm C}^\infty(M)\Rightarrow T\in gl(E)\subset \cl P^1(E,M)
  \end{equation}
  \begin{equation}\label{(v**)}
    \{T\in\cl P(E,M)| [T,\cl P^0(E,M)]\subset \cl P^{i}(E,M)\}=\cl P^{i+1}(E,M).  
  \end{equation}
\end{prop}
\pre
The non-singularity of $\cl P(E,M)$ comes from that of $\cl D(E,M),$ since, $\cl D^1_{\cl A}(E,M)=\cl P^1(E,M).$ 

In addition, $\cl D(E,M)$ being quasi-distinguishing, the relations (\ref{(v*)}) and (\ref{(v**)}) directly come since $\cl P^k(E,M)\subset\cl D^k(E,M),$ for any $k\in\mathbb{N}.$
 
The algebra $\cl P(E,M)$ is symplectic  because $\cl P(E,M)=\cl D(E,M)$. \hfill $\blacksquare$
 
\begin{prop}\label{centralisateur p}
  Let $E\to M$ be a vector bundle of rank $n>1.$ Then for any $\Psi\in\cl C(\cl P(E,M)),$ we have
  \[
    \Psi(\cl P^0(E,M))\subset gl(E)\quad  { and }\quad \Psi(\gamma_u)=\Psi(1)\gamma_u,
  \]
  for any $u\in{\rm C}^\infty(M).$
\end{prop}
\pre
For any $u,v\in{\rm C}^\infty(M), $ we have
 \[
   [\Psi(\gamma_u),\gamma_v]=0.
 \]
Then, in virtue of the relation (\ref{(v*)}) of the Proposition \ref{pe non sing symple pas quasi dist}, $\Psi(\gamma_u)\in gl(E).$ 
In the same way as for the Proposition \ref{centralisateur}, we have for all $ u\in{\rm C}^\infty(M),$
   \[
      \Psi(\gamma_u)=\Psi(1)\gamma_u,
   \]   
because we can come back to a step where the fact that the quantum Poisson algebra $ \cl P (E, M) $ is non-singular allows us to conclude. Indeed, we have the relation
\[
  \Psi(\gamma_u\widehat{T}(\gamma_u))= \gamma_u\Psi(\widehat{T}(\gamma_u)),
\] 
for all $T\in\cl P^1(E,M)$ and $\gamma_u\in\cl P^0(E,M),$ exactly as in the general case of the Proposition \ref{centralisateur}, previously cited. 

Replacing $T$ by $\gamma_v T $ and $\gamma_u$ by $\gamma_u+\gamma_w,$ the equality above induces
 \[
   \Psi((\gamma_u+\gamma_w)\gamma_v \widehat{T}(\gamma_u+\gamma_w))=(\gamma_u+\gamma_w)\Psi(\gamma_v\widehat{T}(\gamma_u+\gamma_w)).
 \] 
Therefore, replacing $\gamma_v$ with $\widehat{T}(\gamma_w),$ a development analogous to that of the Proposition \ref{centralisateur} makes it possible to obtain, for any $u,w\in{\rm C}^\infty(M),$
\[
  \Psi(\gamma_u(\widehat{T}(\gamma_w))^2)=\gamma_u\Psi((\widehat{T}(\gamma_w))^2).
\]
As $ \cl P (E, M) $ is non-singular, we obtain the desired result by applying the same reasoning as for the proposition cited above. \hspace*{1cm}\hfill $\blacksquare$


We propose in the following lines a result of Lie-algebraic characterization of vector bundles with the Lie algebra $\cl P^1(E,M).$ 
 
\begin{theo}\label{P1 induit sur P0}
Let $E\to M, F\to N$ two vector bundles of respective ranks $n,n'>1$ with $H^1(M,\mathbb{Z}/2)=0.$
If $\Phi:\cl P^1(E,M)\to\cl P^1(F,N)$ is an isomorphism of Lie algebras, then
 \[
  \Phi(\cl P^0(E,M))=\Phi(\cl P^0(F,N)).
 \] 
\end{theo}
\pre
The approach consists in proving that $ \cl P^0 (E, M) $ is the Casimir of the subset $Nil(\cl P^1(E,M))$ of $\cl P^1(E,M).$ To achieve this, we show that any element of this Casimir permutes with $sl(E)$ and therefore with $ gl (E) $ whole.
Let us first observe that
 \[
 Nil(\cl P^1(E,M)) \subset  gl(E)
 \]
This comes from the fact that locally, $ T \in \cl P^1 (E, M) $ is written $T=A+\gamma_{u_i}\partial_i$ and  $\cl D(M)$ is distinguishing. 
Indeed, this writing allows to see that for all $ T \in \cl P^1 (E, M), $ relative to a connection on $ E, $ there exists $ X \in Vect (M) $ such that, for any $v\in{\rm C}^\infty(M),$
\[
 (adT)(\gamma_v)=L_X(v)id =\gamma_{X\cdot v}.
\]
We therefore deduce that
\[
 (adT)^p(\gamma_v)= (L_X)^p(v)id,
\]
for any $v\in{\rm C}^\infty(M)$ and any $p\in\mathbb{N}.$ 
The announced result follows since $ \cl D (M) $ is distinguishing.
Let us now calculate the Casimir of $ Nil (\cl P^1(E, M)). $
Recall that we have, by definition, 
 \[
  Cas(Nil(\cl P^1(E,M)))= \{T\in Nil(\cl P^1(E,M)): [Nil(\cl P^1(E,M)),T]=0\}
 \]
We will show the equality
\[
  Cas(Nil(\cl P^1(E,M)))= {\rm C}^\infty(M)id.
\]
 It suffices to prove that, for all $ A\in Cas(Nil(\cl P^1(E,M)))$ and $ B\in sl(E),$ 
\[
   [A_x,B_x]=0 \quad \forall x\in M.
\]
Let $x_0\in M$ and $(U,\psi)$ be a trivialization of the vector bundle $Hom(E,E)\to M,$ whose domain contains $x_0, $ associated with a trivialization $(U,\varphi)$ of $E\to M.$ 
Recall that for all $x\in U,$ $\varphi_x: E_x\to \mathbb{R}^n$ is a linear bijection and that it is the same for 
\[
 \psi_x:gl(E_x)\to gl(n,\mathbb{R}):A\mapsto \varphi_x\circ A\circ\varphi_x^{-1}.
\] 
For $\beta \in gl(E_{x_0}),$ let $B\in gl(E)$ such that 
\begin{equation}\label{eqCasNil}
  \left\{
  \begin{array}{l}
      B_x=0 \mbox{ for } x\in M\setminus U\\
      B_x=u(x)\psi^{-1}_x\circ\psi_{x_0}(\beta) \mbox{ for } x\in U,
  \end{array}
  \right.
\end{equation}
where $u\in {\rm C}^\infty(M)$ is a function with compact support in $ U $ and with value $1$ at $x_0.$
We then have $B_{x_0}=\beta.$ 
Moreover, we observe that if there exists $ k \in \mathbb {N} $ such that $ \beta^k = 0, $ then $ B^k = 0,$ and we can therefore conclude in this case that $B\in Nil(\cl P^1(E,M)).$
Indeed, this comes from the fact that for all $ A, B \in gl (E), $ we have that $ (adB)^k (A) $ decomposes into a sum whose terms are of the form
  \[
   c_{\alpha,\beta}B^\alpha\circ A\circ B^\beta, \quad c_{\alpha,\beta}\in\mathbb{R}, \alpha+\beta=k.
  \]
Now let $ A \in Cas(Nil (\cl P^1 (E, M))) $ and consider a nilpotent element $\beta\in gl(E_{x_0}).$
For an element $B\in gl(E)\cap Nil(\cl P^1(E,M))$ obtained as in (\ref{eqCasNil}) above, we then have  
\[
  0=[A_{x_0},B_{x_0}]=[A_{x_0},\beta].
\]
Which means that $ A_{x_0} $ commutes with all the nilpotent matrices of $ gl (E_{x_0}), $ and therefore with their linear envelope which is $sl(E_{x_0}).$
This being true for all $x_0\in M,$ we obtain the sought result. \hfill $\blacksquare$\\
Contrary to what happens with the entire Lie algebra $ \cl D (E, M) $ (or $\cl P(E,M)$)  it is possible, by virtue of the previous Proposition \ref {P1 induit sur P0}, to obtain a Lie-algebraic characterization of vector bundles without considering $\cl P^1(E,M)$ as a ${\rm C}^\infty(M)-$module.

\begin{theo}\label{iso filtre entre pe et pf}
Let $E\to M, F\to N$ be two vector bundles  of respective ranks $n,n'>1.$ Then any isomorphism $\Phi: \cl P^1(E,M)\to\cl P^1(F,N) $ of Lie algebras is such as its restriction to $\cl P^0(E,M)$ is of the form
 \[
   \Phi|_{\cl P^0(E,M)}=\kappa \Psi,
 \]
where $\Psi:\cl P^0(E,M)\to \cl P^0(E,N)$ is an isomorphism of associative algebras. We also has 
 \[
  \Phi(gl(E))=gl(F).
 \]
\end{theo}
\pre
 In virtue of the previous Proposition \ref{P1 induit sur P0}, we have
\[
 \Phi(\cl P^0(E,M))=\cl P^0(F,N).
\]
We also have, for any element $A\in gl(E),$  
\[
\Phi\circ\gamma_A\circ\Phi^{-1}\in\cl C(\cl P^1(F,N)).
\]
In the above relation, notations are those specified in the proof of the Proposition \ref{iso modules alg quasi filtré}. 
And then for any $\gamma_w\in\cl P^0(F,N),$ the following equality is true
 \begin{equation}\label{(*x)}
   \Phi\circ\gamma_A\circ\Phi^{-1}(\gamma_w)=\Phi(A\Phi^{-1}(1))\gamma_w 
 \end{equation}
The previous Proposition \ref{centralisateur p}  then gives 
\[
 \Phi(A\Phi^{-1}(1))\in gl(F)
\]
for any $A\in gl(E).$
Therefore, $\cl P^1(E,M)$ being symplectic, we deduce that
 \[
   \Phi(gl(E))=gl(F),
 \]
because then, the constant $ \Phi^{-1} (1) $ is not zero.
The rest of the proof is an adaptation of the reasoning developed for the Proposition \ref{iso modules alg quasi filtré}.
Indeed, for $A=\gamma_u,$ $u\in{\rm C}^\infty(M)$ and $\lambda:=\Phi^{-1}(1),$ the  previous relation (\ref{(*x)})  gives in particular 
\[
  \Phi(\gamma_u\Phi^{-1}(\gamma_w))=\Phi(\lambda\gamma_u)\gamma_w\cdot
\]
Setting, $\Phi^{-1}(\gamma_w)=\gamma_v,$ we obtain
 \[
   \Phi(\gamma_u\cdot\gamma_v)=\lambda\Phi(\gamma_u)\Phi(\gamma_v). 
 \]
Therefore, the map
 \[
   \gamma_u\mapsto \Phi^{-1}(1)\Phi(\gamma_u)
 \] 
 is effectively an isomorphism of $\mathbb{R}-$algebras of $\cl P^0(E,M)$ to $\cl P^0(F,N).$ \hfill $\blacksquare$\\

We can now state the following result.

\begin{theo}\label{theo pe caract fibré}
Let $E\to M, F\to M$ be two vector bundles of respective ranks  $n,n'>1$ with $H^1(M,\mathbb{Z}/2)=0.$
The Lie algebras  $\cl P^1(E,M)$ and $\cl P^1(F,M)$  are isomorphic if and only if the vector bundles $E$ and $F$ are. 
\end{theo}

For the Lie algebra $\cl P(E,M)$ we can now define the space of symbols as in  \cite{GraPon1,GraPon2,GraPon4}, for instance.


\section{The classical Poisson algebra $\cl S(\cl P(E,M))$}
 
The  classical limit of the quantum Poisson algebra $ \cl P (E, M)$ is defined by
  \[
     \cl S(\cl P(E,M))=\oplus_{i\in\mathbb{Z}}\cl S^i(\cl P(E,M)),
  \]  
with $\cl S^i(\cl P(E,M))=\cl P^i(E,M)/\cl P^{i-1}(E,M).$

We obtain a "classical Poisson algebra", see \cite{GraPon1}, whose operations are given in the following.
Let us recall that for any $T\in\cl P^i(E,M),$ $ord(T)=i,$ if $T\notin \cl P^{i-1}(E,M).$
For $i\geq ord(T),$ the $i$-degree symbol of $ T $  is defined by
\[
\sigma_i(T)=\left\{
\begin{array}{l}
 0 \mbox{ if } i>ord(T)\\
 T+\cl P^{i-1}(E,M)\mbox{ if } i=ord(T)
\end{array}
\right.
\]
The symbol related to the quantum Poisson structure of $ \cl P (E, M)$ is for its part given by
  \[
  \sigma_{pson} : \cl P(E,M)\to  \cl S(\cl P(E,M)): T\mapsto \sigma_{ord}(T).
  \]
For $P\in \cl S^i(\cl P(E,M))$ and $Q\in\cl S^j(\cl P(E,M))$ such that $P=\sigma_i(T)$ and $Q=\sigma_j(D),$ we set, by definition,
 \[
  P.Q=\sigma_{i+j}(T\circ D) \mbox{ and } \{P,Q\}=\sigma_{i+j-1}([T,D])\cdot
 \]
 
\subsection{Particular case of  the Lie sub-algebra $gl(E)\subset\cl P^1(E,M)$}\label{cas part sle}
 
By virtue of the calculations made in the previous section, we have, by definition, 
 \[
 \sigma_{pson}(\gamma_u)=\gamma_u+\{0\},\quad \forall u\in{\rm C}^\infty(M)
 \]
and we will simply denote $\qquad\sigma_{pson}(\gamma_u)=\gamma_u.$\\

Likewise, for $A\in gl(E)\setminus\cl P^0(E,M),$ we have
\[
 \sigma_{pson}(A)= A'+\cl P^0(E,M),
\] 
with $A'=A-\frac{tr(A)}{n}id,$\quad $tr(A)$ being the trace of $A.$ 
Therefore, for any $A,B\in gl(E),$ we have 
\begin{eqnarray*}
 \sigma_{pson}([A,B]) & = & [A,B]+ \cl P^0(E,M)\\
                      & = & [A',B']+\cl P^0(E,M).
\end{eqnarray*} 
And for the product, if $A,B\notin\cl P^0(E,M),$ we have 
  \[
  \sigma_{pson}(A)\cdot\sigma_{pson}(B)=0.
  \]
For $\gamma_u\in\cl P^0(E,M),$ we then have   
  \[
  \sigma_{pson}(\gamma_u)\cdot\sigma_{pson}(A)=\gamma_u \circ A'+\cl P^0(E,M).
  \]
We therefore have the following identification of Lie algebras
\[
\sigma_{pson}(gl(E))\cong sl(E)\oplus {\rm C}^\infty(M)id,
\]
where the multiplication is commutative and defined by
 \[
 (A+\gamma_u)\cdot(B+\gamma_v)=\gamma_v\circ A +\gamma_u\circ B +\gamma_{uv}
 \] 
and the bracket being given by the following relation
\[
  \{A+\gamma_u , B+\gamma_v \}=[A,B].
\]

\subsection{General case}

Let us state the following result which gives the local expression of the elements of $ \cl P (E, M) $ in a trivialization of $ E. $

Sometimes, for convenience of writing, we simply denote $ \gamma_u $ by $ u \in {\rm C}^\infty (M).$
\begin{prop}\label{critere pkde }
The elements of $\cl P^k(E,M),$ $k\geq1,$ are characterized by the fact that they are written locally, in a trivialization domain $ U \subset M, $ in the form
\begin{equation}\label{(*v)}
  \sum_{|\alpha|<k}T_\alpha\partial^\alpha+\sum_{|\beta|=k}u_\beta\partial^\beta
\end{equation}
where $T_\alpha\in{\rm C}^\infty (U,gl(n,\mathbb{R}))$ and $u_\beta\in{\rm C}^\infty (U).$
\end{prop}
\pre
Let $ T \in \cl P^1(E, M) $ and suppose that in $ U $ we can write
  \[
    T=A+\sum_{1\leqslant i\leqslant m}T_i\partial_i,\quad A,T_i\in{\rm C}^\infty (U,gl(n,\mathbb{R}))\cdot
  \]    
The relation $[T,\gamma_u]\in\cl P^0(E,M)$ gives 
 \[
  \sum_{1\leqslant i\leqslant m} T_i\circ\partial_i(u)\in{\rm C}^\infty(U),\quad\forall u\in{\rm C}^\infty(M).
 \]
Therefore, for any $i\in[1,m],$ we have  $T_i\in{\rm C}^\infty(U).$ 
Assume the result is true for any element of $\cl P^r(E,M)$, with $r<k,$ and let $T\in\cl P^k(E,M)$ such that we locally have
 \[
 T=\sum_{|\alpha|<k }T_\alpha\partial^\alpha+\sum_{|\beta|=k}T_\beta\partial^\beta  .
\] 
We have, by applying the induction hypothesis to $ [T, \gamma_u], $ that its highest order term of derivation is of the form
 $
 \sum_{|\lambda|=k}u_\lambda\partial^\lambda.
 $
Now in $ [T_\beta \partial^\beta, \gamma_u], $ the highest order terms of derivation are of the form $\partial_iu T_{\beta}\partial^{\beta_i},$ with $\beta_i=\beta-e_i,$ where $(e_i)_{1\leq i\leq m}
$ designates the canonical basis of the $\mathbb{R}-$vector space $\mathbb{R}^m.$ 
We therefore have, for any $u,$ $\partial_{i}uT_{\beta}\in{\rm C}^\infty(U),$ and consequently $T_{\beta}\in{\rm C}^\infty(U).$
Conversely if $ T \in\cl D (E, M) $ is written locally in the form (\ref{(*v)}), we obtain that $T\in\cl P^k(E,M),$ for all $k\geq 1,$ by doing a recurrence on $k.$
More precisely, a such element $T$  has the following form 
\begin{equation}\label{(**)}
 \sum_{|\alpha|=k-1}A_{\alpha}\partial^\alpha+\sum_{|\beta|=k}u_{\beta}\partial^\beta.
\end{equation}
with $A_\alpha\in {\rm C}^\infty(U,sl(n,\mathbb{R}))$  and $u_\beta\in{\rm C}^\infty(U).$
\hfill $\blacksquare$\\


We notice that the local decomposition (\ref{(**)})  given above is not intrinsic.
In fact, if in the sum on the right we recognize the principal symbol of the differential operator in the usual sense, the sum on the left, for its part, does not resist a change of coordinates and is therefore not globally defined.

In the following lines we build a global decomposition allowing to find a global meaning to the expression given in (\ref{(**)})  previously.\\

Let now $T\in\cl S^{k-1}(M)\otimes sl(E)$ and assume given a unit partition $(U_i,\rho_i)$ of $M$ whose domains $U_i$ are the trivialisation one of $E.$ In any $U_i,$ if $T$ is expressed in the form
 \[
   T=\sum_{|\alpha|=k-1}A_{\alpha,i}\xi^\alpha\cdot
 \] 
We then set
 \[
  \overline{T}_i=\sum_{|\alpha|=k-1}A_{\alpha,i}\partial^\alpha\in\cl D^{k-1}(E|_{U_i},U_i)
 \]
with $A_{\alpha,i}\in  {\rm  C}^\infty(U_i,sl(n,\mathbb{R})).$ 
The differential operator
  \[
    \overline{T}=\sum_i \rho_i\overline{T}_i\in \cl D^{k-1}(E,M)\subset\cl P^k(E,M),
  \]
associated with the partition of the unit chosen at the start is then such that
   \[
    \sigma_{pson}(\overline{T})=\sigma_{ppal}(\overline{T})=T,
   \]
where $\sigma_{ppal}$ is the usual principal symbol. Note that $\overline{T}$ is not unique.\\
Nevertheless, we have the following statement.
\begin{prop}\label{prop suite deg k}
The space $\cl S^k(\cl P(E,M))=\cl P^k(E,M)/\cl P^{k-1}(E,M)$ of symbols in the sense " quantum Poisson algebra " of differential operators in $\cl P^k(E,M)$ is determined by the following  exact short sequence of $\mathbb{R}-$vector spaces 
   \begin{equation}\label{suite deg k}
    0\longrightarrow\cl S^{k-1}(M)\otimes sl(E)\stackrel{\theta}{\longrightarrow}\cl P^k(E,M)/\cl P^{k-1}(E,M)\stackrel{\delta}{\longrightarrow}\cl S^k(M)\longrightarrow 0,
   \end{equation}
with $\cl S(M)=Pol(T^*M)$, 
$\theta: T\mapsto \overline{T}+\cl P^{k-1}(E,M)$ and 
  \[
  \delta: D+\cl P^{k-1}(E,M)\mapsto \left\{
                                 \begin{array}{l}
                                 0 \mbox{ if } D\in\cl D^{k-1}(E,M)\\
                                 \sigma_{ppal}(D) \mbox{ if not. }
                                 \end{array}
                                \right. 
  \] 
\end{prop}
\pre
The application $\theta$ is well-defined. Indeed, the differential operators $D_1,D_2\in\cl D^{k-1}(E,M)\subset\cl P^k(E,M)$ are such that $\sigma_{pson}(D_1)=\sigma_{pson}(D_2),$ we then do have $D_1-D_2\in\cl P^{k-1}(E,M);$ which means that the image of $ T $ does not depend on the choice of the operator $\overline{T}$ such that $\sigma_{pson}(\overline{T})=T.$ \\ 
Also, $ \theta $ is obviously a linear map and it is injective. Indeed, let $T\in\cl  S^{k-1}(M)\otimes sl(E)$ such that $\theta(T)=0.$ 
We then have $\overline{T}\in\cl P^{k-1}(E,M).$ 
But by construction $\overline{T}\notin\cl D^{k-2}(E,M).$ We thus have
    \[
      \sigma_{ppal}(\overline{T})\in\cl S^{k-1}(M)id
    \] 
We deduce, since $\sigma_{pson}(\overline{T})=\sigma_{ppal}(\overline{T})=T,$ that $T=0.$ 
 The map $ \delta $ being linear and directly surjective, we show to finish that
\[
ker(\delta)=Im(\theta).
\] 
The inclusion $ker(\delta)\supset Im(\theta)$ is obvious. 
Let's prove the other sense of that inclusion. If $D+\cl P^{k-1}(E,M)\in ker(\delta),$ then, by the definition of $\delta,$ we have
 \[
   D\in \cl D^{k-1}(E,M)\cap\cl P^k(E,M).
 \]	
Consequently, we obtain
 \[
  \sigma_{ppal}(D)\in \cl S^{k-1}(M)\otimes gl(E)\cap \cl S^{k-1}(M)\otimes sl(E);
 \]  
and the inclusion sought is a direct result of this.\hfill $\blacksquare$\\
Therefore, seen as $\mathbb{R}-$vector spaces, we have the following decomposition 
  \[
    \cl S^k(\cl P(E,M))=Pol^{k-1}(T^*M,sl(E))\oplus Pol^k(T^*M,\mathbb{R})
  \]
for any integer   $k\in\mathbb{N}.$ But for the Lie algebra structure, the above decomposition is not true because the following exact sequence of Lie algebras (but also of  associative algebras) is not split. 
\begin{equation}\label{(v)}
 \xymatrix{
 \quad 0\ar[r] & \cl S(M)\otimes sl(E)\ar[r] & \cl S(\cl P(E,M))\ar[r] & \cl S(M)\ar[r]& 0 
         }
\end{equation}
But for first order differential operators, we have the following result.  
\begin{prop}\label{suite scindee}
Let $E\to M$ be vector bundle of rank $n.$  
With respect to a connection on $ E, $ the following short exact sequence of Lie algebras is split
\[
\xymatrix{       
  0\ar[r] & \cl P^0(E,M)\stackrel{i}{\longrightarrow} \cl P^1(E,M)\stackrel{\sigma_{pson}}{\longrightarrow} \cl P^1(E,M)/\cl P^0(E,M)\ar[r] & 0 
         }
\]
where $i$ is the canonical injection and $\sigma_{pson}$ the map previously defined.
\end{prop}
\pre
Via a covariant derivation $ \nabla $ of $ E, $ we have the identification of $ \mathbb {R}-$vector spaces 
  \begin{equation}\label{idP1}
    \cl P^1(E,M)\simeq Vect(M)\oplus gl(E).
  \end{equation}
Indeed, for $T\in\cl P^1(E,M)\setminus sl(E),$ we have 
\[
  \sigma_{ppal}(T)\simeq  X\in Vect(M),
\] 
with the identification $\cl S^1(M)id\simeq Vect(M).$ Therefore, since $\nabla_X\in\cl P^1(E,M),$ the difference $T-\nabla_X$ is an endomorphisms field. 
We denote by  $\lambda$ the linear bijection given in (\ref{idP1}).

For $ T=\nabla_X+A$ and $D=\nabla_Y+B,$ we then get in $\cl P^1(E,M)$
  \[
   [T,D]=\nabla_{[X,Y]}+R^\nabla(X,Y)+ \nabla_XB-\nabla_YA+[A,B]\cdot
  \] 
This is equivalent to
 \begin{equation}\label{(*)}
  \quad [(X,A),(Y,B)]=([X,Y],R^\nabla(X,Y)+ \nabla_XB-\nabla_YA+[A,B])
 \end{equation}
in the space $Vect(M)\oplus gl(E).$   
Moreover, consider the following short exact sequence
\[
\xymatrix{       
  0\ar[r] & sl(E)\ar[r] & \cl P^1(E,M)/\cl P^0(E,M)\ar[r] & Vect(M)\ar[r] &  0 
         }
\]
It is corresponding to the particular case $ k = 1 $ of (\ref{suite deg k}).
Therefore, as in Proposition \ref{prop suite deg k}, we have the map $\delta$ 
   \[
  \delta: T+\cl P^{0}(E)\mapsto \left\{
                                 \begin{array}{l}
                                 0 \mbox{ if } T\in gl(E)\\
                                 \sigma_{ppal}(T) \mbox{ if not }
                                 \end{array}
                                \right. 
  \] 
which is surjective and the injection $\theta: A\in sl(E)\mapsto A+\cl P^0(E,M).$
Seen as $\mathbb{R}-$vector spaces, we then have the following identification
 \[
    \cl P^1(E,M)/\cl P^0(E,M)\cong Vect(M)\oplus sl(E). 
 \]
Consider the following commutative diagram
\[
  \xymatrix{(X,A)\in Vect(M)\oplus sl(E) \ar[d]^{\lambda^{-1}}\ar[r]^\mu & \cl P^1(E,M)/\cl P^0(E,M)\\  
  \nabla_X+A\in\cl P^1(E,M)\ar[ur]_{\sigma_{pson}} &
           }
\]
The linear map $ \mu = \sigma_{pson} \circ \lambda^{- 1} $ is injective because $ \mu (X, A) = 0 $ induces $ \nabla_X + A \in \cl P^0 (E); $ and we deduce that
\[
X=0 \quad\mbox{ and } \,A\in sl(E)\cap\cl P^0(E,M).
\]   
Likewise, $ \mu $ is surjective.
The bracket in $\cl P^1(E,M)/\cl P^0(E,M)$ is given by
\[
   [[T],[D]]  = \nabla_{[X,Y]}+R^\nabla(X,Y)+ \nabla_XB-\nabla_YA+[A,B]+\cl P^0(E,M)
 \]
with $[T]=\nabla_X+A+\cl P^0(E,M)$ and $[D]=\nabla_Y+B+\cl P^0(E,M).$\\
Therefore, the corresponding operation in $ Vect (M) \oplus sl (E), $ obtained by structure transport via $ \mu, $ is not necessarily a Lie bracket since the term $ R^\nabla ( X, Y) $ is not always of zero trace. \\
To remedy this, suppose that $ \nabla $ is associated with a connection form of a reduction of the frame principal bundle $ L^1(E) $ of $ E $ to the Lie subgroup $ O (n) $ of $ GL (n, \mathbb {R}). $ For such a derivation, $ R^\nabla $ has values in $ sl (E). $\\
We then have an isomorphism of Lie algebras
\[
\mu : Vect(M)\oplus sl(E)\to \cl P^1(E,M)/\cl P^0(E,M),
\]
the space $ Vect (M) \oplus sl (E) $ being provided with the following bracket
 \begin{equation}\label{(***)}
   [(X,A),(Y,A)]=([X,Y],R^\nabla(X,Y)+ \nabla_XB-\nabla_YA+[A,B]).
 \end{equation}  
Given the relations  (\ref{(*)})  and  (\ref{(***)}), we conclude that the canonical injection
\[
 \beta: Vect(M)\oplus sl(E)\to Vect(M)\oplus gl(E)
\]
is a homomorphism of Lie algebras.
Therefore, let us consider the map
\[
  \lambda^{-1}\circ \beta\circ \mu^{-1}:\cl P^1(E,M)/\cl P^0(E,M)\to \cl P^1(E,M)\cdot
\]
It is a homomorphism of Lie algebras allowing to identify the Lie algebra $ \cl P^1 (E, M) / \cl P^0 (E, M) $ to a Lie subalgebra of $ \cl P^1 (E, M), $ for the structures specified in the previous lines, and we can see that it is a section of $ \sigma_{pson}. $ \\
We have just shown that the short exact sequence of the statement is split. \\
\hspace*{4cm}\hfill $\blacksquare$

Note that with the notations in the previous proof, for $ T = \nabla_X + A, $ the decomposition
\[
T=(\nabla_X+A-\frac{1}{n}tr(A)\,id)+\frac{1}{n}tr(A)\,id\,
\]
only depends on the reduction and not on the connection choice.
Indeed, if relatively to another covariant derivation $ \nabla', $ associated with the same reduction, we consider an analogous decomposition of $ T, $ then we have
 \begin{eqnarray*}
   T=\nabla'_X+A'& = & \nabla_X+(A'+(\nabla'_X-\nabla_X))\\
                 & = & \nabla_X+S+A' 
 \end{eqnarray*}
with $S=\nabla_X-\nabla'_X,$ and thus the trace of $S$ is null.

We conclude that $A'=A-S$. Therefore, $tr(A)=tr(A')$ and we have
  \[
   \nabla'_X+A'-\frac{1}{n}tr(A')\,id=\nabla_X+S+A-S-\frac{1}{n}tr(A)\,id.
  \] 



 

\section{Lie-algebraic characterization of vector bundles}

\begin{theo}
Let $E\to M, F\to M$ be vector bundles of respective ranks $n,n'>1$ with $H^1(M,\mathbb{Z}/2)=0.$
The Lie algebras $ \cl S (\cl P (E, M)) $ and $ \cl S (\cl P (F, M)), $ seen as $ {\rm C}^\infty (M) -$modules, are isomorphic if and only if the vector bundles $ E $ and $ F $ are.
\end{theo}
\pre
Let $\Phi : \cl S (\cl P (E, M))\to \cl S (\cl P (F,M)) $ be an isomorphism of Lie algebras.
Note that $\Phi$ preserves $\cl A={\rm C}^\infty(M)$ the basis of these classical Poisson algebras.
Indeed, for any $u\in {\rm C}^\infty(M),$ we have
\[
  \Phi(\gamma_u)=u\Phi(\gamma_1).
\]
And the conclusion comes from the fact that the quantum Poisson algebras considered are symplectic.
Let now $B\in gl(E).$ We then have
$
  \{B,\cl A \}=0.
$
Therefore $\{\Phi(B),\cl A\}=0$ and this implies $\Phi(B)\in gl(E),$ since $\cl D(E,M)$ is quasi-distinguishing, according to the Proposition \ref{alg quasi P q non-sing,sympl, quasi-q}. 
\hfill $\blacksquare$\\

For the Lie subalgebras $ \cl S^1 (\cl P (E, M)) $ and $ \cl S^1 (\cl P (F, N)), $ this result may improve. To prove this, we use the short exact sequence presented in the previous Proposition \ref{suite scindee}.

\begin{theo}
Let $E\to M, F\to M$ be two vector bundles of respective ranks $n,n'>1$ with $H^1(M,\mathbb{Z}/2)=0.$
The Lie algebras $ \cl S (\cl P^1(E, M)) $ and $ \cl S (\cl P^1(F, M)), $  are isomorphic if, and only if, the vector bundles $ E $ and $ F $ are. 
\end{theo}
\pre
We observe the decomposition
\[
\cl S^1(\cl P(E,M))=sl(E)\oplus Vect(M)
\] 
obtained previously via a connection on $ E.$ 
For all $ T = (\nabla_X, A),$ let $ B \in sl (E) $ such that there exists $ r \in \mathbb {N} $ verifying 
\[
 (ad(T))^r(B)=0.
\]
We then have, if $\nabla_X$ is $r$ times applied,
 \[
   \nabla_X(\nabla_X\cdots(\nabla_X(B)))=0\cdot
 \]
 In a trivialization of $E$ of domain $U\subset M,$ considering $B$ whose local expression is of the form $(\alpha_{ij})=(\delta_{12}u),$ $u\in {\rm C}^\infty(U),$ i.e. having all its terms null except that of the position $(1,2),$ we can choose $u$ so that we necessarily have $X=0.$ 
We deduce that
  \[
   Nil(\cl S^1(\cl P(E))\subset sl(E).
  \] 
Moreover, we know that if $A\in sl(E)$ is such that $A^p=0,$ we then have $A\in Nil(\cl S^1(P(E,M)).$ 
Indeed, observe that
\begin{eqnarray*}
 (ad A)^k(\nabla_X+B) & = & (ad A)^k(\nabla_X)+(ad A)^k(B) \\
                      & = & -(ad A)^{k-1}(X\cdot A)+ (ad A)^k(B)
\end{eqnarray*} 
where, by virtue of the particular case studied in the previous section \ref{cas part sle}, for all $C,D\in sl(E),$ we have
\begin{eqnarray*}
  (ad\, C)(D) & = & \{C,D\} =  [C,D].       
\end{eqnarray*}
This allows us to conclude, since $ (ad \, C)^k (D) $ is then a sum of the terms of the form $a_k C^k\circ D\circ C^{k-1}, a_k\in\mathbb{R}$.
In the following lines, the goal is to establish that
\[
\left.  \right\rangle  Nil(\cl S^1(\cl P(E,M))\left\langle= sl(E), 
  \right.
\]
where the usual notation $\left.  \right\rangle  H\left\langle\right.$ designates the linear envelope of the subset $H$ of a vector space.
Let $ A \in sl (E). $ Over a trivialization domain $ U \subset M, $ we can therefore write
  \[
    A|_U=\Sigma_{i}N_i^U  
  \]
with $N_i^U \in sl(E|_U),  (1\leq i\leq n^2-1)$, which are nilpotent endomorphism fields since $sl(n, \mathbb {R}) $ admits a basis formed of nilpotent matrices.
Consider now a Palais cover of $M,$ 
\[
  \cl O=\cl O_1\cup\cdots\cup\cl O_r, r\in \mathbb{N},
\]
locally finite, the elements $ U _{\alpha,j} $ of each $ \cl O_j $ being trivialization domains of $ E $ 2 by 2 disjoint, and  $(\cl \rho_{\alpha,j}),$ a partition of the unit, locally finite and subordinate to this cover. 
We therefore have 
\begin{eqnarray*}
  \rho_{\alpha,j}A & = & \sum_{i=1}^{n^2-1} N_{i,\alpha,j}
\end{eqnarray*} 
with $N_{i,\alpha,j}\in sl(E|_{U_{\alpha,j}})$ nilpotent and compactly supported in $U_{\alpha,j}$. \\
Let now pose $\cup\cl O_j=\cup_{\alpha}U_{\alpha,j}=U_j$ and consider $N_{ji}$ defined by
\[
     N_{ji}(x)=  \left\{ \begin{array}{l}
           N_{i,\alpha,j}(x) \mbox{ if } x\in U_{j}\\
           0 \quad\quad\quad \mbox{ if not }
               \end{array} 
               \right.
\] 
We then obtain the smoothness of $N_{ji}$. Indeed, for any $x\notin U_j,$ consider an open neighborhood $V\ni x$ of compact adherence. We know that $ V $ is only encountered by a finite number of supports of $N_{\alpha,j},$ whose reunion is the compact that we agree to denote by $K.$ 

Therefore $V\setminus K$ is an open neighborhood of $ x $ in which $ N_{ji} $ is identically zero. And we have thus established our assertion since for $ x \in U_j $ the smoothness of  $N_{ji}$ is obvious. 
We conclude that 
\[
  \sum_{\alpha} \rho_{\alpha,j} A=\sum_{i=1}^{n^2-1} N_{ji}
\mbox{ and }  
   A=\sum_{i=1}^{n^2-1}\sum_{j=1}^rN_{ji}.
 \]
We deduce that any isomorphism 
 $
 \Phi: \cl S(\cl P^1(E,M))\to  \cl S(\cl P^1(F,N))
 $
of Lie algebras is necessarily such that
\[
  \Phi(sl(E))= sl(F).
\]
Hence, by virtue of the  Theorem \ref{ende carct e}, we have the desired result.
\hspace*{1cm}\hfill $\blacksquare$

\end{document}